\documentclass[12pt,a4paper]{article} 
\usepackage{url}
\usepackage{graphicx}
\usepackage{amsmath}

\newcommand{\C}{\mbox{$I \!\!\!\! C$}}

\newcommand{\Ei}{\operatorname{Ei}}
\newcommand{\Ci}{\operatorname{Ci}}
\newcommand{\Si}{\operatorname{Si}}
\newcommand{\si}{\operatorname{si}}

\newcommand\PVint{\mathop{\setbox0\hbox{$\displaystyle\intop$}%
    \hskip0.2\wd0%
    \vcenter{\hrule width0.6\wd0height0.5pt depth0.5pt}%
    \hskip-0.8\wd0%
    }\mskip-\thinmuskip\intop\nolimits}

\title{On a Singular Integrodifferential Equation arising from a Linearised Free Surface Problem}
\author{Patrick McLean\thanks{School of Mathematica and Physics, University of Tasmania, Private Bag 37, Hobart, Tasmania, 7001, \textsc{Australia}. \protect\url{mailto:p_mclean@maths.utas.edu.au}}}

\date{(31 November 2005)}

\begin{document}

\maketitle

\begin{abstract}
A problem of linear surface waves discussed by Forbes \cite{forbes2} initially gave rise to a singular integrodifferential equation over the real line. We have been able to transform this integrodifferential equation into a linear second order differential equation whose solution has been found explicitly in terms of the sine and cosine integral functions. Furthermore, we have been able to recover known results of physical interest using more simple techniques.
\end{abstract}

\tableofcontents

\section{Introduction}
In this paper we consider the singular integrodifferential equation
\begin{equation}
\label{eq:ide}
u(x)+F^2Hu'(x) =f(x), \qquad -\infty<x<\infty,
\end{equation}
where \( H \) is the Hilbert transform, \( F > 0  \) is a constant and \( f(x) \) is a certain algebraic function. We show that ~(\ref{eq:ide}) is equivalent to a linear second order differential equation which can be explicitly solved in terms of the sine and cosine integral functions.

This equation arises from a  free surface problem in fluid mechanics and this background is discussed in section 2. In section 3 we use certain properties of the Hilbert transform to convert ~(\ref{eq:ide}) to a differential equation the solution of which is then evaluated in terms of the sine and cosine integral functions. 
We are also able to extract certain properties of the surface. Finally, in section 4 we present graphs of the solution \( u(x) \) for a variety of \( F \).

\section{Background}
In \cite{forbes1} and \cite{forbes2} Forbes considers surface waves caused by flow over a (line) vortex. He models the situation using a complex potential \( \Phi(z) \) beneath a free surface \( y=S(x) \). This results in a nonlinear free surface problem. Linearisation (assuming the vortex strength is small) and employing the boundary integral method results in the following singular integrodifferential equation
\begin{equation}
\label{eq:ide0}
\frac{\partial\phi_1}{\partial x} (x,0) + \frac{F^2}{\pi} \PVint_{-\infty}^\infty \frac{\partial^2\phi_1}{\partial^2 t} (t,0) \frac{dt}{t-x} = \frac{1}{\pi(x^2+1)}, \quad -\infty<x<\infty,
\end{equation}
where \( \phi_1 \) is the real part of the linearised complex potential and the integral is to be interpreted in the Cauchy principal value sense. This equation appears as equation (2.8) in \cite{forbes2}.

If we define the Hilbert transform of a function \( f(x) \) as in \cite[\S 15]{erdelyi}
\begin{equation}
Hf(x) = \frac{1}{\pi} \PVint_{-\infty}^\infty \frac{f(y)}{y-x} dy,
\end{equation}
then we can express equation ~(\ref{eq:ide0}) in the form ~(\ref{eq:ide}) where
\begin{equation}
u(x)= \frac{\partial\phi_1}{\partial x}(x,0), \qquad f(x) = \frac{1}{\pi(x^2+1)}.
\end{equation}
The linearised surface and complex potential can be recovered from \( u(x) \).
 
Equation ~(\ref{eq:ide}) is an integrodifferential equation, thus to solve it uniquely we need to specify an extra condition. Namely, we seek a solution to equation ~(\ref{eq:ide}) which vanishes towards \( -\infty \), that is, we shall impose the boundary condition
\begin{equation}
\label{eq:ic}
\lim_{x\rightarrow-\infty}u(x)=0.
\end{equation}

\section{Analytical Solution}
In this section we convert the integrodifferential equation ~(\ref{eq:ide}) to a differential equation from which we obtain various results regarding the solution \( u(x) \). In particular, we obtain a representation of the function \( u \) in terms of sine and cosine integral functions,  we determine the asymptotic behaviour of \( u(x) \) towards \( +\infty \) and we evaluate \( u(0) \) and \( u'(0) \).

\subsection{Conversion to a differential equation\label{sec:tode}}
We can convert the integrodifferential equation ~(\ref{eq:ide}) into a second order differential equation. The resulting equation can be solved and the boundary condition ~(\ref{eq:ic}) invoked to give a representation of \( u(x) \) as an integral. 

Differentiating ~(\ref{eq:ide}), applying \( H \) and multiplying by \( F^2 \) gives
\begin{equation}
\label{eq:ide1}
F^2 Hu' + F^4 H(Hu')' = F^2 Hf'.
\end{equation}
Using the results \( H^2u=-u \), see \cite[\S 15.1 (2)]{erdelyi} and \( Hu'=(Hu)' \), see \cite[\S 15.1 (8)]{erdelyi} we have that \( H(Hu')' = -u'' \). Thus, ~(\ref{eq:ide1}) becomes
\begin{equation}
  \label{eq:ide2}
  F^2 Hu' - F^4 u'' = F^2 Hf'. 
\end{equation}
Subtracting ~(\ref{eq:ide2}) from ~(\ref{eq:ide}) gives the 2nd order differential equation
\begin{equation}
  \label{eq:de}
  F^4 u'' + u = f-F^2 Hf'.
\end{equation}

The general solution to the homogeneous part of this equation is
\begin{equation}
A \sin(\frac{x}{F^2}) +B \cos(\frac{x}{F^2})
\end{equation}
for arbitrary constants \( A \) and \( B \). By the method of variation of parameters, see, for example, \cite{boycediprima}, we find that a particular solution, \( u_p \) say, which is zero at \( -\infty \) is given by
\begin{equation}
u_p(x) = \frac{1}{F^2} \int_{-\infty}^x \sin(\frac{x-t}{F^2}) \left[f(t)-F^2 Hf'(t) \right] dt.
\end{equation}
Thus, the general solution of equation ~(\ref{eq:de}) is given by
\begin{equation}
u(x)=A \sin(\frac{x}{F^2}) +B \cos(\frac{x}{F^2})+u_p(x).
\end{equation}
Imposing condition ~(\ref{eq:ic}) we have, since \( u_p(-\infty) =0 \), that both \( A \) and \( B \) must be zero. So, the required solution is 
\begin{equation}
\label{eq:varpar}
u(x) =\frac{1}{F^2} \int_{-\infty}^x \sin(\frac{x-t}{F^2}) \left[f(t)-F^2 Hf'(t) \right] dt.
\end{equation}
Since \( Hf' = (Hf)' \) we can apply integration by parts to give 
\begin{equation}
\label{eq:twopar}
\int_{-\infty}^x \sin(\frac{x-t}{F^2}) (Hf)'(t) dt = \left[ \sin(\frac{x-t}{F^2}) Hf(t) \right]_{-\infty}^x + \frac{1}{F^2} \int_{-\infty}^x \cos(\frac{x-t}{F^2}) Hf(t) dt.
\end{equation}
The following Hilbert transform is to be found in \cite[\S 15.2 (10)]{erdelyi}
\begin{equation}
\label{eq:htransalg}
H(\frac{1}{x^2+1}) = \frac{-x}{x^2+1},
\end{equation}
from which we see that
\begin{equation}
\label{eq:vanish}
\lim_{t\rightarrow -\infty} Hf(t) = 0.
\end{equation}
From ~(\ref{eq:varpar}), ~(\ref{eq:twopar}), and  ~(\ref{eq:vanish}) we have the representation
\begin{equation}
\label{eq:sol2}
u(x) = \frac{1}{F^2} \int_{-\infty}^x \left[ \sin(\frac{x-t}{F^2}) f(t) + \cos(\frac{x-t}{F^2}) Hf(t) \right] dt.
\end{equation}

\subsection{Trigonometric Integrals}
We follow the standard reference work \cite[\S 5.2]{as} and define the trigonometric integrals as follows.
For \( z \in \C \) define the sine integral \( \Si(z) \) as
\begin{equation}
\Si(z) = \int_0^z \frac{\sin t}{t} dt.
\end{equation}
The sine integral \( \Si \) is an entire function. A commonly used notation is \(\si(z) = \frac{\pi}{2} -\Si(z) \).

For \( z \in \C \) such that \( | \arg z | < \pi \) define the cosine integral \( \Ci(z) \) as
\begin{equation}
\Ci(z) = \gamma + \log z + \int_0^z \frac{\cos t-1}{t} dt.
\end{equation}
The cosine integral \( \Ci \) has a branch cut discontinuity along the negative real axis.

The functions \( \Si \) and \( \Ci \) occur as Fourier sine and cosine transforms.
Specifically, for \( |\arg a| < \pi \) and \( y>0 \) we have that
\begin{eqnarray}
\label{eq:cosint}
\int_0^\infty \frac{\cos(xy)}{a+x} \, dx & = & -\si(ay)\sin(ay) -\Ci(ay)\cos(ay) \\ 
\label{eq:sinint}
\int_0^\infty \frac{\sin(xy)}{a+x} \, dx & = & \Ci(a y)\sin(ay) -\si(a y)\cos(ay),
\end{eqnarray}
see \cite[\S 1.1 (9)]{erdelyi} and \cite[\S 1.2 (10)]{erdelyi}, respectively.

\subsection{Evaluation of Solution in Terms of Trigonometric Integrals}
In this section we evaluate \( u(x) \) in terms of the trigonometric integrals \( \Si \) and \( \Ci \). 

Invoking the change of variable \( x-t=s \) in equation ~(\ref{eq:sol2}) we obtain
\begin{equation}
u(x) =\frac{1}{\pi F^2} \int_0^\infty \left[\sin(\frac{s}{F^2}) \frac{1}{(x-s)^2+1}+\cos(\frac{s}{F^2})\frac{x-s}{(x-s)^2+1} \right] ds.
\label{final}
\end{equation}

On making partial fraction expansions of the algebraic expressions in equation ~(\ref{final}) and using the identities ~(\ref{eq:cosint}),~(\ref{eq:sinint}), together with standard trigonometric identities we obtain the representation
\begin{eqnarray}
u(x) & = & \frac{\sin(\frac{x}{F^2}) }{\pi F^2e^{\frac{1}{F^2}} } 
\left( \frac{1}{2i} \left[ \Ci(-\frac{x-i}{F^2})-\Ci(-\frac{x+i}{F^2}) \right] -\frac{1}{2} \left[\Si(-\frac{x-i}{F^2}) +\Si(-\frac{x+i}{F^2}) \right] + \frac{\pi}{2} \right) \nonumber \\
& & + \frac{\cos(\frac{x}{F^2}) }{\pi F^2e^{\frac{1}{F^2}} } \left( \frac{1}{2} \left[ \Ci(-\frac{x-i}{F^2})+\Ci(-\frac{x+i}{F^2}) \right] +\frac{1}{2i} \left[\Si(-\frac{x-i}{F^2}) -\Si(-\frac{x+i}{F^2}) \right] \right).
\end{eqnarray}
This representation allows evaluation of \( u(x) \) for arbitrary \( x \) in terms of the functions \( \Si \) and \( \Ci \) at complex arguments. High precision evaluation of these functions is provided by computer packages such as Mathematica. For Fortran code for evaluating the sine and cosine integrals at complex arguments see \cite{acm}. 

An alternative method of evaluating \( u(x) \) is to view ~(\ref{final}) as a linear combination of the Fourier sine and cosine transforms
\begin{equation}
\int_0^\infty \sin(\frac{s}{F^2}) \frac{1}{(x-s)^2+1} ds \qquad
\int_0^\infty \cos(\frac{s}{F^2})\frac{x-s}{(x-s)^2+1}ds,
\end{equation}
and use the method of \cite{paper2} to numerically evaluate them.

\subsection{Asymptotic Behaviour of Solution}
In this section we examine the asymptotic behaviour at \( +\infty \) of the solution ~(\ref{eq:sol2}).
Rewrite ~(\ref{eq:sol2}) as
\begin{eqnarray}
u(x)&=& \frac{1}{F^2} \int_{-\infty}^\infty \left[ \sin(\frac{x-t}{F^2}) f(t) -\cos(\frac{x-t}{F^2}) Hf(t) \right] dt \\
 & &-\frac{1}{F^2} \int_{x}^{\infty} \left[ \sin(\frac{x-t}{F^2}) f(t) -\cos(\frac{x-t}{F^2}) Hf(t) \right] dt.
\end{eqnarray}
The second integral in this representation vanishes as \(x\rightarrow\infty \), so we have that
\begin{equation}
u(x) \sim \frac{1}{F^2} \int_{-\infty}^\infty \sin(\frac{x-t}{F^2}) f(t) dt - \frac{1}{F^2} \int_{-\infty}^\infty \cos(\frac{x-t}{F^2}) Hf(t) dt,
\label{uatinf}
\end{equation}
for large positive \( x \). We shall need the following two properties of the Hilbert transform:
\begin{equation}
\int_{-\infty}^\infty f(x-y) Hg(y) dy  = \int_{-\infty}^\infty Hf(x-y) g(y) dy.
\label{eq:hconv}
\end{equation}
and 
\begin{equation}
\label{eq:htranstrig}
H\cos(x) = -\sin(x),
\end{equation}
which are to be found in \cite[\S 4.2 (16)]{tricomi} and \cite[\S 15.2 (47)]{erdelyi}, respectively.

Thus, ~(\ref{uatinf}) becomes
\begin{eqnarray}
u(x) & \sim & \frac{1}{F^2} \int_{-\infty}^\infty \sin(\frac{x-t}{F^2}) f(t) dt - \frac{1}{F^2} \int_{-\infty}^\infty H\cos(\frac{x-t}{F^2}) f(t) dt \\
& = & \frac{2}{F^2} \int_{-\infty}^\infty \sin(\frac{x-t}{F^2}) f(t) dt.
\end{eqnarray}
Now, using the fact that the sine function is odd and \( f \) is even it follows that
\begin{equation}
u(x) \sim \frac{2}{\pi F^2} \sin(\frac{x}{F^2}) \int_{-\infty}^\infty \cos(\frac{t}{F^2}) \frac{1}{t^2+1} dt.
\end{equation}
On using the Fourier cosine transform to be found in \cite[\S 1.2 (11)]{erdelyi}
\begin{equation}
\int_0^\infty \frac{\cos(xy)} {x^2+a^2} \, dt = \frac{\pi}{2 a} e^{-ay}, \qquad a,y>0
\end{equation}
we obtain
\begin{equation}
u(x) \sim 2 e^{-\frac{1}{F^2}} \sin(\frac{x}{F^2}), \quad x \rightarrow \infty.
\label{eq:asypmt}
\end{equation}
This agrees with \cite[eqn 3.4]{forbes1}. Thus, we have found the asymptotic behaviour at \( +\infty \) of the function \( u(x) \). That is, the solution oscillates sinusoidally at \( +\infty \).

\subsection{Behaviour at origin}
We now note the value of \( u(x) \) and its derivative at the origin \( x=0 \). Evaluating ~(\ref{final}) at \( x=0 \) we have
\begin{equation}
u(0) = \frac{1}{\pi F^2} \int_0^\infty \frac{\sin(\frac{s}{F^2})-s\cos(\frac{s}{F^2}) }{s^2+1}ds.
\end{equation}
The following identities appear in \cite{erdelyi} as equations 2.2 (14) and 1.2 (12) respectively
\begin{eqnarray}
\int_0^\infty \frac{\sin(xy)}{x^2+a^2} \, dx & = & \frac{1}{2a} \left[ e^{-ay} \Ei(ay) - e^{ay} \Ei(-ay) \right] \\
\int_0^\infty \frac{x\cos(xy)}{x^2+a^2} \, dx & = & -\frac{1}{2a} \left[ e^{-ay} \Ei(ay) + e^{ay} \Ei(-ay) \right],
\end{eqnarray}
where \( a, y >0 \) and the exponential integral function \( \Ei \) is defined as
\begin{equation}
\label{eq:expintei}
\Ei(z) = -\PVint_{-z}^\infty \frac{e^{-t}}{t} \, dt,
\end{equation}
see \cite[eqn 5.1.2]{as}. Thus, we are able to express \( u(0) \) as
\begin{equation}
u(0) = \frac{e^{-\frac{1}{F^2}}}{\pi F^2} \Ei(\frac{1}{F^2}).
\end{equation}

Similarly, we can derive that
\begin{equation}
u'(0)  =  \frac{e^{-\frac{1}{F^2}}}{F^4}.
\end{equation}

\section{Linearised Surface Profile}
In this section we present graphical representations of the surface profile for three Froude numbers \( F \). From equations (2.1) and (2.5a) of \cite{forbes1} the surface profile is given by
\begin{equation}
S(x)=-\epsilon F^2 \frac{\partial \phi_1}{\partial x}(x,0)+ O(\epsilon^2).
\end{equation}

\begin{figure}[tbp]
\includegraphics*{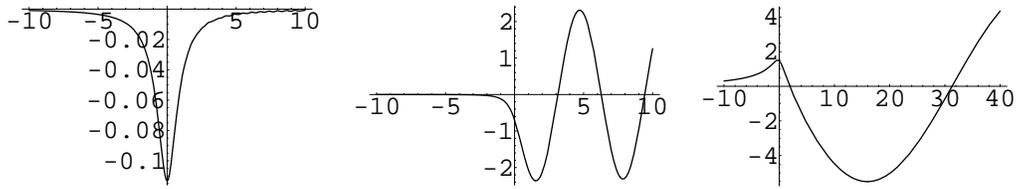}
\caption{Surface profiles for \(F=0.1,1,10\) }
\label{waves}
\end{figure}

We take \( \epsilon =1 \), \( F=0.1,1,10 \) and plot \( u(x) \) in figure \ref{waves}. We note that in each case \( u(x) \) is a superposition of an impulse due to the vortex and a sinusoidal wave downstream. Qualitatively the main differences are scale and the relative strengths of these two contributions.

\end{document}